\documentclass[12pt,english,reqno]{smfart}
\usepackage[T1]{fontenc}

\usepackage[english,francais]{babel}

\usepackage{amssymb,url,xspace,smfthm}
\usepackage{bm,upgreek}
\usepackage{extarrows} 
\usepackage{enumitem}
\usepackage{url}
\usepackage{graphicx}

\title[]{Discrete splicing theorem for noise sensitivity of invasion percolation}
\date {6 May 2013 (retyped 20 October 2018)}

\author[K.N.S.]{Kamron Nirou Saniee}
\address{New York, NY}
\email{saniikami@gmail.com} 


\begin{document}

\begin{abstract}
We state and prove a version for invasion percolation of Schramm and Smirnov's (2011) discrete gluing theorem for critical planar percolation. The result is a first step toward establishing factorization and characterizing noise sensitivity for invasion.
\end{abstract}

\maketitle

\tableofcontents

\section{Motivation: noise sensitivity, black noise}
\noindent There has been considerable recent interest in noise sensitivity of spatial random processes, and of percolation models in particular. Noise sensitivity of an event describes the property that, by resampling a portion of the state configuration of the process, the realization of the event is rendered nearly independent of its original status. One may also say that with high probability, a configuration containing any small level of random error gives nearly no predictive measure of the occurrence of the event. The current formulation \cite{Ga11,EF12} of this property stems from two distinct directions of inquiry: that of the robustness in the above sense of Boolean functions arising in computer science, which motivated use of the discrete Fourier spectrum of events (\cite{GPS10}) to characterize sensitivity of percolation events, and the study initiated by Tsirelson and Vershik \cite{TV98, Ts04} of a special class of stochastic process which arise from continuous products of probability spaces, and which in some sense possess intrinsically the property of sensitivity. \\

\noindent Tsirelson's theory of noises provides a framework for characterizing certain stochastic processes which generalize sequences of i.i.d random signs. The underlying  process (the \textit{noise}) is translation invariant and independent on disjoint domains. One is often interested in functionals (or observables) that capture global properties of the ambient process, such as macroscopic events in statistical physics models. Important examples of spatial noises include critical planar percolation (established in \cite{SS11}) and the Brownian web, and in 1D, white noise (or the derivative of Brownian motion) and Watanabe's coalescent flow \cite{W01}. We refer to \cite{Ts} for the comprehensive construction of this theory. \\

\noindent As remarked in \cite{SS11}, noise provides one definition of the scaling limit of critical planar percolation. By arguments of \cite{Ts}, characterizing the ambient process (such as percolation) as a noise provides one route to characterizing functionals (e.g., the crossing event) as sensitive or stable. We emphasize that this approach identifies stability or sensitivity as a byproduct of the ambient process rather than as property of the specific functional, as indicated in the extreme case of \textit{black noises} for which all functionals are by default sensitive. The main difficulty in proving a process is a noise is establishing that it \textit{factorizes} on adjacent domains, in the sense that the sigma algebras generated by these domains determine the sigma-algebra of the union. Once factorization is shown, sensitivity can follow from bounds on the probability of events which control the functional of interest. This is made precise in the below overview of basic results on noises, and critical exponents known prior in the physics literature provide such bounds for percolation \cite{SW01}. In Section 3 we state the essential results in Schramm and Smirnov's proof of factorization and in Section 4 we prove discrete splicing, a mesh-dependent analogue on annular factors, for the invasion percolation model.

\section{Noise and stable $\sigma$-field}
\noindent Here we recall some facts about the noise theory framework from \cite{Ts}. \\

\begin{defi} A \textit{continuous factorization} of a probability space $(\Omega,\mathcal{F},\mathbb{P})$ is a collection of sub-sigma-fields $\mathcal{F}_D\subset \mathcal{F}$ for all $d$-dimensional rectangles $D\subset\mathbb{R}^d$ such that
\begin{center}
$\mathcal{F}_D\otimes\mathcal{F}_{D'}=\mathcal{F}_{D''}$
\end{center} whenever $D,D'$ form a partition of $D''$, and such that
the collection $(\mathcal{F}_D)_D$ generates $\mathcal{F}$. \\
\end{defi}

\begin{defi}
A \textit{noise} consists of a probability space $(\Omega,\mathcal{F},\mathbb{P})$, a continuous factorization $(\mathcal{F}_D)_D$ of $\mathcal{F}$ and a measurable action $(T_h)_{h\in\mathbb{R}^d}$ of the additive group of $\mathbb{R}^d$ on $\Omega$, such that $\mathcal{F}_D\underset{T_h}{\mapsto}\mathcal{F}_{h+D}$ for all $h\in \mathbb{R}^d$. \\ 
\end{defi}

\begin{exem}[Brownian noise] Consider $k$ random signs $\tau_k$ assigned to points $k\epsilon$, $-M\leq k\leq M$. The associated probability space is
\begin{center}
$(\Omega_{\epsilon,M}=\{+1,-1\}^{\epsilon\mathbb{Z}\cap[-M,M]},\mathcal{F}_{\epsilon,M},P_{\epsilon,M})$
\end{center}
where $P_{\epsilon,M}(A)=\frac{|A|}{|\Omega_{\epsilon,M}|}$ and $\mathcal{F}_{\epsilon,M}$ is the sigma-field generated by the signs. Taking a Riemann Integrable function $\phi\in L^2(\mathbb{R})$ with $\int_{\mathbb{R}} \phi^2 dx=1$, the sum
\begin{center}
$\epsilon^{1/2}\sum_k \phi(k\epsilon)\tau(k\epsilon)$
\end{center}
converges in distribution (for $\epsilon\rightarrow 0$, $M=1/\epsilon$) to a $N(0,1)$ random variable which we note $\int_{\mathbb{R}} \phi(x) dB(x)$. Brownian motion $B(\cdot)$ follows as
\begin{center}
$B(x)=\int_\mathbb{R} \mathbf{1}_{[0,x]} dB(x)$ if $x>0$, 0 if $x=0$, $-\int_\mathbb{R} \mathbf{1}_{[x,0]}dB(x)$ if $x<0$. 
\end{center} 

\noindent The limiting probability space is described by $\Omega=\mathbb{R}^{\infty}$ (the space of increments of Brownian motion), $P=N(0,1)^{\otimes\infty}$ where $P(\{(\alpha_1,\alpha_2,\ldots):\alpha_1<a_1,\ldots, \alpha_d<a_d\})=N(0,1)((-\infty,a_1])\cdots N(0,1)((-\infty,a_d])$ for $N(0,1)$ the gaussian law, and $\mathcal{F}$ the sigma-field of $N(0,1)^{\otimes\infty}$-measurable sets (i.e. the union of Borel subsets and measure-zero sets). 
Equivalently, $\mathcal{F}$ is the sigma field generated by the random variables $B(y)-B(x)$ for $x<y$. Let $\mathcal{F}_{a,b}$ be the sigma-field generated by $B(y)-B(x)$ for $(x,y)\subset (a,b)$. By the defining properties of Brownian motion, \\

\begin{enumerate}[label=(\roman*)]
\item $W_0=0$
\item $W_t$ is almost-surely continuous, and 
\item increments on disjoint intervals are independent and distributed $W_t-W_s \sim N(0, t-s)$, \\
\end{enumerate}

\noindent a continuous factorization is given by
\begin{equation*}
\mathcal{F}_{a,c}=\mathcal{F}_{a,b}\otimes\mathcal{F}_{b,c}.
\end{equation*}

\noindent Formally, the Brownian noise (or white noise) is identified as the isometric map
\begin{equation*}
L^2(\mathbb{R})\ni \phi \mapsto \int \phi dB \in L^2(\Omega,\mathcal{F}, P).\\ \medskip
\end{equation*}
\end{exem}

\noindent\textbf{Stability  } $L^2(\Omega,\mathcal{F},P)$ certainly contains random variables other than the linear $\int_{\mathbb{R}} \phi(x) dB(x)$ (that is, limits of linear functions of the signs $\tau(k\epsilon)$). 
Let $\omega_\epsilon\in\Omega_{\epsilon,M}$ be the configuration $(\tau(k\epsilon))_{-M/\epsilon\leq k\leq M/\epsilon}$ and let $H_n^{\epsilon}$ be the space of sums of the form 
\begin{equation*} 
X_{\epsilon,n}=\epsilon^{n/2}\sum_{-M/\epsilon\leq k_1<\ldots<k_n\leq M/\epsilon} \psi(k_1\epsilon,\ldots,k_n\epsilon)\tau(k_1\epsilon)\cdots\tau(k_n\epsilon)
\end{equation*}
and denote by $H^n$ the space of limits $(\lim_{M\rightarrow\infty}\lim_{\epsilon\rightarrow 0}(\ldots))$ of these sums.
For $A\subset [-M/\epsilon,M/\epsilon]$, define the random variable
\begin{equation*}
N_pX_{\epsilon,n}\overset{\Delta}{=}X_{\epsilon,n}(N_p\omega)
\end{equation*}
where $N_p\omega$ is obtained from $\omega$ by resampling each sign $\tau_k$ for $k\in A$ with probability $p$ (i.e., by independently changing signs in $A$ with probability $p$). \\

\noindent The variance $||X_{\epsilon,n}-N_p X_{\epsilon,n}||^2$ is a measure of the sensitivity of $X_{\epsilon,n}$ to the local resampling. This motivates the following definition: \\

\begin{defi}
We call \textit{stable sigma-field}, or $\mathcal{F}_{stable}$, the sigma-field generated by $\{X\in L^2(\Omega,\mathcal{F},P): ||X-N_p X||\xrightarrow[p\to 0]{}0\}$. \\
\end{defi}
 
\noindent An equivalent interpretation is that information about the underlying noise process coming from the resampled region is corrupted by some level $p$ of noise. Stability and sensitivity refer to the qualitative property of a functional $X$ that one may reliably observe its realization in the presence of noisy data. Sensitivity is equivalent to the decorrelation, for arbitrarily small levels of noise, of the true and noisy observable. \\

\begin{prop}
Let $\sigma(X)=\mathcal{F}_{s,t}$. Then $X\in H_1$ if, for all $u\in(s,t)$, 
\begin{equation*}
X=\mathbb{E}(X|\mathcal{F}_{s,u})+\mathbb{E}(X|\mathcal{F}_{u,t})
\end{equation*}

\end{prop}
\begin{prop}
We have $\sigma(H_1)=\mathcal{F}_{stable}$.
\end{prop}

\begin{coro}

Let $L^n=\{s_0^n,s_1^n,\ldots,s_n^n\}\subset\mathbb{R}$ be subsets of $(a,b)$ such that $s_0^n<\ldots<s_n^n$ for each $n$ and such that $\cup_{n\geq 1} L^n$ is dense in $(a,b)$. Then the orthogonal projection \\
\begin{equation*}
\left \{ X\in L^2(\Omega,\mathcal{F},P):\mathbb{E}(X)=0 \right \} \hookrightarrow \left \{ X\in H_1: X\: \text{ }\mathrm{is}\text{ }\:\mathcal{F}_{a,b}\mathrm{-measurable} \right \} \\
\end{equation*}
\\
is given by
\begin{equation}
\mathrm{Proj}_{H_1}=\lim_{n\rightarrow\infty}\; \sum_{k=1}^{n} \mathbb{E}(\cdot |\mathcal{F}_{s_{k-1}^n,s_k^n}) 
\end{equation}
\end{coro}

\begin{defi}
A random variable $X\in L^2(\Omega,\mathcal{F},P)$ is \textit{sensitive} if
\begin{equation*}
\mathrm{Proj}_{H_1}X=0.
\end{equation*}
\end{defi}

\noindent It follows that the space of observables $L^2(\Omega,\mathcal{F},P)$ splits into two orthogonal subspaces $L^2(\mathcal{F}_{\text{stable}})$ and $L^2(\mathcal{F}\setminus \mathcal{F}_{\text{stable}})$. Orthogonality may be seen by noting that, for $n\neq m$, $X_{n,\epsilon}\in H_n^{\epsilon},X_{m,\epsilon}\in H_m^{\epsilon}$,
\begin{equation*}
\mathrm{Cov}(X_{n,\epsilon},X_{m,\epsilon})=0
\end{equation*}
by independence of the signs $\tau(k\epsilon)$. \\

\noindent The following corollary to (1) provides a useful criterion for sensitivity. \\

\begin{theo}[{[BKS99]}]
Let $(\mathcal{F}_{x,y})_{x\leq y}$ be a continuous factorization over $\mathbb{R}$ and $X\in L^2(\Omega,\mathcal{F},P)$. Then $X$ is orthogonal to $L^2(\mathcal{F}_{stable})$ if
\begin{equation*}
\underset{\{s_1,\ldots,s_n\} \uparrow}{\limsup} \sum_{k=1}^{n+1}\mathrm{Var}(\mathbb{E}(X|\mathcal{F}_{s_{k-1},s_k}))
\end{equation*}
where the $\limsup$ is taken over all finite subsets of $\mathbb{R}$ ordered by inclusion (i.e., sets $L^n$ in the above projection theorem, with $(a,b)=(-\infty,+\infty)$).\\
\end{theo}

\noindent A noise for which all elements of $L^2(\Omega,\mathcal{F},P)$ are orthogonal to $H_1$ is referred to as \textit{black noise}. All $L^2$ functionals of a black noise are sensitive. \\

\begin{coro}
Let $D\in \mathbb{R}^2$ be a rectangle and $x\in D^{\circ}$. Consider a two-dimensional noise (i.e., with factorization over the algebra of planar rectangles) such that for all $X\in L^2(\Omega,\mathcal{F},P),$
\begin{equation}
\mathbb{E}(X|\mathcal{F}_{x+[-\epsilon,\epsilon]^2})^2=o(\epsilon).
\end{equation}
Then the noise $((\Omega,\mathcal{F},P),(\mathcal{F}_D)_{D\subset\mathbb{R}^2})$ is black. \\
\end{coro} 

\begin{exem}[Brownian web]
A \textit{Brownian web} on a probability space $(\Omega,\mathcal{F},P)$ is a measurable map
\begin{equation*}
\phi: \Omega\times \{(s,t)\in\mathbb{R}^2,s\leq t\}\times \mathbb{R} \rightarrow \mathbb{R}
\end{equation*}
\begin{equation*}
(\omega,(s,t),x)\mapsto \phi_{(s,t)} (t)(\omega)
\end{equation*} 
such that for every finite collection $(s_i,x_i)$ ('starting points' of the web), the processes $\phi_{(s_i,x_i)}(\cdot)$ ('trajectories' of the web) form a system of coalescing Brownian motions: that is, a set
\begin{equation*}
(X^1,\ldots,X^n)
\end{equation*}
 such that the $X^i$ start from $x_i$ at time $s_i$ and are independent until the earliest time $T$ where $\; i\neq j$ such that $X^i(T)=X^j(T)$; for $t\geq T$, $X^i(t)=X^j(t)$ and $(X^k)_{k\leq i}$ form a system of $n-1$ coalescing Brownian motions. \\
 
\noindent Let $\mathcal{F}_{x,y}$ be the sigma-field generated by the restriction of trajectories of the web to the vertical strip $(x,y)\times (-\infty,\infty)$. It follows from the definition of Brownian motion that $(\mathcal{F}_{x,y})_{x\leq y}$ is a continuous factorization (the vertical factorization of the web). Ellis and Feldheim \cite{EF12} show the association of sigma-fields to strips $(-\infty,\infty)\times (x,y)$ gives a continuous factorization (horizontal factorization). In \cite{Ts}, Tsirelson shows (2) for $((\Omega,\mathcal{F},P),(\mathcal{F}_{x,y})_{x\leq y})$. Together these results show that the Brownian web factorized on planar rectangles is a 2-dimensional black noise. \\
\end{exem}

\begin{rema}
The Brownian web and critical percolation are the only known examples of two-dimensional black noise. \\
\end{rema}

\begin{rema}
The quantity $\mathbb{E}(X|\mathcal{F}_{x+[-\epsilon,\epsilon]^2})^2$ is referred to as the \textit{influence} of the cell $x+[-\epsilon,\epsilon]^2$. Bounding the influence of small resampled regions provides one route to showing blackness. This is the method employed for the proof for the Brownian Web and is implicit in Smirnov's proof for critical percolation. In the case of binary or integer-valued observables (such as the percolation crossing event), a cell is said to be \textit{pivotal} if the noise, restricted to this cell, controls an event whose occurrence changes the value of the observable. For both of these models, exact critical exponents for pivotal events do not need to be known in advance to characterize the underlying process as a noise, and the main difficulty lies in proving that one may reconstruct observables from adjacent spatial domains. \\
\end{rema}

\begin{rema}
For a given model, the boundary of the partition of $D$ into disjoint subsets is subject to regularity requirements to ensure that the sigma-field of the union may be recovered from the smaller sigma-fields. For example, it is claimed in \cite{SS11} that critical site percolation on the triangular lattice factorizes as a noise for boundary curves of Hausdorff dimension at most $5/4$. Stronger conditions on factors for the Brownian Web is an open problem. Below, we make the same regularity assumption as \cite{SS11} on the boundary curve.
\end{rema}

\section {Invasion Percolation Model}

\subsection{Invasion and critical percolation}\;\; Invasion percolation has received recent attention as system exhibiting "self-organizing criticality": it possesses no endogenous parameter, yet defines an object which closely resembles the infinite critical percolation cluster. Statistical, connective and geometric properties of the invaded cluster coincide with the critical cluster  \cite{Ja03,DSV09,Zh95}. The two models are however structurally distinct, as invasion does not share translational invariance. While it may be defined on a general graph, we restrict to invasion percolation on the lattice $\mathbb{Z}^2$ in this presentation and in our main result.\\

\begin{defi}[Invasion Percolation Cluster.]
\noindent Let $(\tau_e)_e$ be independent random variables uniformly distributed in $[0,1]$, indexed by edges (or bonds) $e\in E^2=E(\mathbb{Z}^2)$ of the lattice $\mathbb{Z}^2$. We refer to the $\tau_e$ as \textit{weights}. Define the boundary of a subgraph $G\subset (\mathbb{Z}^2,E^2)$ as 
\begin{equation*}
\Delta G=\left \{e=(x,y)\in E: e\notin E(G) \wedge (x\in G \vee y\in G) \right \}.
\end{equation*}
Now let $G_0=(\{0\},\varnothing)$ and $(G_i)_{i=1}^{\infty}$ be a sequence of subgraphs of the lattice induced by the edge sets 
\begin{equation*}
E(G_{i+1})= E(G_i)\cup\underset{e\in \Delta G_i}{\text{argmin}} \left \{\tau_e \right \}.
\end{equation*}
The \textit{invasion percolation cluster} (IPC) is defined as $\cup_{i=0}^{\infty} \; G_i$.\\ \medskip
\end{defi}

\noindent Hence, the invasion graph starts growing from the origin, and adds the lowest-weight boundary edge at each step. $G_i$ is referred to as the invasion at time $i$. We denote by $\Omega\subset\{0,1\}^{E^2}$ the set of configurations accessible by this procedure, $\mathcal{F}=\sigma(\tau_e: e\in E^2)$ the sigma-algebra generated by the weights, and $\mathbb{P}$ the induced measure on $\Omega$. We write $\omega\sim\mathbb{P}$ to indicate a sample $\omega$ from $\mathbb{P}$.\\

\noindent The above construction provides a coupling of the invasion process to standard planar percolation with parameter $p$, by assigning open (closed) status to $e$ if $\tau_e\underset{(\geq)}{<}p$. The distribution of the resulting graph of $p$-open edges coincides with that of the usual definition, where each edge is open (closed) with probability $p$ (resp. $1-p$), independently of other edges. Hence, in what follows, invaded bonds are referred to also as open bonds. We write $A\underset{{}^{\mathbf{\bullet}}}{\leftrightarrow} B$ to indicate two sets $A$ and $B$ are connected by edges having the property "$	\bullet$" (e.g., "\textit{p-open}"). We will also refer to configurations on the \textit{dual lattice} $(\mathbb{Z}^2)^*$,

\begin{equation*}
(\mathbb{Z}^2)^*=(1/2,1/2)+\mathbb{Z}^2
\end{equation*}
as follows: for $\omega\in \Omega$, define the \textit{dual configuration} $\omega^*$ 
\begin{equation*}
\omega^*(e^*)=\tau_e
\end{equation*}
by assigning the weight $\tau_e$ to the edge in $E(\mathbb{Z}^2)^*$ which intersects $e$. A dual edge $e^*$ is $p$-open if $e$ is $p$-open. \\

\subsection{Factorization of critical planar percolation} Let 
\begin{equation*}
(\Omega_p=\{\text{0,1}\}^{\eta E^2},\mathcal{F}_p,\mathbb{P}_p=(p\delta_1+(1-p)\delta_0)^{\otimes \eta E^2})
\end{equation*}

\noindent be the probability space of standard percolation with parameter $p$ in the mesh-$\eta$ square lattice, $\eta\mathbb{Z}^2$, and $\omega\sim \mathbb{P}_p$. Let $Q$ be a \textit{quad} (a homeomorphic image of $[0,1]^2$) in $\mathbb{R}^2$ and designate two opposite sides of $Q$. A quad is \textit{crossed by }$\omega$ if the restriction of $\omega$ to edges intersecting $Q$ contains an open cluster which intersects the designated sides. A theorem of Kesten \cite{Ke80} implies that as the mesh $\eta$ of the lattice tends to zero, for fixed $Q$, the probability of a crossing tends to 0 (resp. 1) if $p\underset{(>)}{<}$the \textit{critical parameter }$p_c=1/2$. This is one formulation of the phase transition between the subcritical and supercritical regimes $p\underset{(>)}{<}1/2$, which originates in the experimental heuristic of looking for crossings of macroscopic regions in order to tell whether the system is subcritical or supercritical. By a result of \cite{CCN85}, for any $p>p_c$, the IPC almost surely intersects the infinite $p$-open cluster. The above definition implies that once the invasion enters the cluster, it does not leave. Combined with nonexistence of the infinite cluster for $p<1/2$, it follows that if $e_i$ denotes the $i^{th}$ invaded edge, then $\underset{i\rightarrow\infty}{\limsup}\;\tau_i=p_c$. For further relations between the IPC and percolation clusters, we refer to \cite{DSV09}. \\

\noindent In \cite{SS11}, Smirnov and Schramm show that critical planar percolation factorizes as a noise, for factors belonging to the algebra of planar domains generated by rectangles. As planar percolation is translation invariant, factorization suffices to show it is a noise in the sense of Tsirelson. This reduces to showing that, in a rectangle $R$ with a smooth path $t$ cutting it, for any $\epsilon>0$, there is a finite number (depending on $\epsilon$) of percolation crossing events measurable with respect to $\mathcal{F}_{R\setminus t}$ such that, conditional on these events, the crossing status of $R$ may be predicted within error $<\epsilon$. By \textit{smooth} is meant a set having finite one-dimensional upper Minkowski content $m^{*1}(t)=\underset{\epsilon\rightarrow 0^+}{\limsup} \; \text{area}\{z: \text{dist}(z,t)<\epsilon\}$, a condition which will be motivated below. \\

\begin{theo}[Mesh-independent sampling, {[SS11 4.1]}] 
Let $\mathbb{P}_{\eta}$ be the measure of critical percolation in the mesh-$\eta$ lattice, $Q$ be a piecewise-smooth quad in $\mathbb{R}^2$, and $t$ be a set with $m^{*1}(t)<\infty$. Denote by $\boxminus_Q$ the indicator of the crossing event for $Q$. Then for every $\epsilon>0$, there is a finite set of piecewise smooth quads $Q_\epsilon\subset t^c$ and a set of crossing events $\mathcal{W}_\epsilon$ measurable with respect to $\sigma(\boxminus_{Q'}:Q'\in Q_\epsilon)$, such that
\begin{equation*}
\underset{|\eta|\rightarrow 0}{\lim}\;\mathbb{P}_{\eta}(\mathcal{W}_\epsilon\Delta\boxminus_Q)<\epsilon.
\end{equation*}
\end{theo}

\noindent The proof relies on a discrete version of the result, which states that, uniformly in the mesh size, the crossing status of $Q$ in the discrete percolation configuration may be reconstructed with high probability by sampling crossing events outside a neighborhood of $t$. \\

\begin{theo}[Discrete gluing, {[SS11 1.5]}]
Let $Q$ and $t$ be as above, and let $\mathcal{F}_{\epsilon}$ be the sigma field generated by the restriction of the realization $\omega\sim\mathbb{P}_{\eta}$ to the complement in $\mathbb{R}^2$ of the $\epsilon$-neighborhood of $t$. Then for every $\delta>0$,
\begin{equation*}
\underset{\epsilon \searrow 0}{\lim}\; \underset{\eta\in (0,\epsilon)}{\sup}\; \mathbb{P}_{\eta}(\delta<\mathbb{P}_{\eta}(\boxminus_{Q}|\mathcal{F}_{\epsilon})<1-\delta)=0.
\end{equation*}
\end{theo}

\noindent 

\noindent The purpose of this note is to state and prove a version of discrete gluing for invasion percolation. We consider the functional analogue to $\boxminus_Q$ for the invasion process, namely, the maximal number of disjoint paths connecting the inner and outer boundaries of an annulus $\mathsf{Ann}(n/2,n)$. For a subset $T\subset E^2$ of edges and $\omega\sim\mathbb{P}$, define a random configuration $\omega'$ by 
\begin{equation*}\begin{cases}
\omega'(T)\sim \mathbb{P}(\cdot |\mathcal{F}_{E^2\setminus E(T)}) \\
\omega'(e) =\omega(e) \text{ for any edge } e\notin T.
\end{cases}\end{equation*}\\

\noindent This \textit{resampled configuration}
changes the status of edges in $T$ according to the invasion measure conditional on events outside of $T$. Other resampling procedures may be considered, for example, resampling edge weights $\tau_e$ for $e\in E(T)$ uniformly from $[0,1]$ in a deterministic or random fashion (e.g., as in dynamical percolation, with edges changing weight according to an exponential Poisson clock \cite{SS10}). In this case, the IPC is grown again each time an edge changes. However, the given resampling procedure for $\omega'$ corresponds to the \textit{noising} of the path $t$ implicit in the mesh-independent and discrete gluing results (i.e., by conditioning on events away from the path $t$), and is also the natural one insofar as we are interested in how local connectivity data influences macroscopic properties of the invasion. \\

\noindent Our main result states that if we resample the status of edges contained in a thin annular tranche in the interior of an annulus, then the discrete invasion's behavior outside the tranche accurately predicts the number of disjoint crossings of the annulus. The annular tranche is chosen for convenience, as our argument applies more generally for the $\epsilon$-neighborhood of any closed path $t$ in $S(n)$ enclosing $S(n/2)$ with the regularity condition that it may be covered by at most $\text{const}(t)/\epsilon$ balls of radius $\epsilon$. This is the case if $m^{*1}(t)<\infty$ and provides a bound on the number of resampled events which the argument will need. \\

\section{Main theorem}
\noindent In this section we state and prove the discrete splicing theorem for invasion percolation. \\

\begin{theo}[Discrete splicing for invasion]
Let $\omega\sim\mathbb{P}$ be a sample from the invasion measure and let $\boxtimes_M(\omega)$ denote the event that $\omega$ has a maximal number $M$ of disjoint paths connecting $\partial S(n/2)$ to $\partial S(n)$ in $\omega$. Let $T_\epsilon=\partial S(3n/4)+[-\epsilon/2,\epsilon/2]$ be a thin annular strip about the origin, and $\mathcal{F}_{\epsilon}$ the sigma-field generated by weights $\tau_e$ for edges $e$ intersecting $T_{\epsilon}^c$. Then for any $\delta>0$,
\begin{equation}
\underset{\epsilon \searrow 0}{\lim}\; \underset{n>1}{\sup}\;\mathbb{P}(
\delta<\mathbb{P}(\boxtimes_M|\mathcal{F_{\epsilon}})<1-\delta)=0.
\end{equation}

\end{theo}
 
\noindent\textbf{Proof Strategy } Let $N=N_n(\omega)$ be the maximal number of disjoint paths connecting $\partial S(n/2)$ to $\partial S(n)$ in $\omega$ and $\omega_{\epsilon}$ be a resampling of the configuration $\omega$ in $T_{\epsilon}$. We show
\begin{equation}
\underset{n>1}{\sup\;}\mathbb{P}(N(\omega)\neq N(\omega_\epsilon))\xrightarrow[\epsilon\rightarrow 0]{} 0
\end{equation}

\noindent which implies 
\begin{equation*}
\underset{n>1}{\sup\;}\mathbb{P}(\omega\in \boxtimes_M,\omega_{\epsilon}\notin \boxtimes_M)\xrightarrow[\epsilon\rightarrow 0]{} 0.
\end{equation*} Since $\omega$ and $\omega_{\epsilon}$ are independent conditional on $\mathcal{F}_{\epsilon}$, \begin{equation*}
\mathbb{P}(\omega\in \boxtimes_M,\omega_{\epsilon}\notin \boxtimes_M)=\mathbb{E}[\mathbb{P}(\boxtimes_M|\mathcal{F_{\epsilon}})(1-\mathbb{P}(\boxtimes_M|\mathcal{F_{\epsilon}}))].
\end{equation*} 
By Markov's inequality,
\begin{equation*}
\mathbb{P}(\delta<\mathbb{P}(\boxtimes_M|\mathcal{F_{\epsilon}})<1-\delta)\leq \frac{\mathbb{E}\mathbb{P}(\boxtimes_M|\mathcal{F_{\epsilon}})(1-\mathbb{P}(\boxtimes_M|\mathcal{F_{\epsilon}}))}{\delta^2}
\end{equation*}
and the theorem statement follows.\\

\noindent \textbf{Lemmas}\;\; For $p>p_c$, define the \textit{finite scaling length} $L(p,\delta)$ as
\begin{equation*}
\underset{n\geq 0}{\min}\; \{ \mathbb{P}(\{0\}\times[0,n]\underset{p-\text{open path} \text{ in } [0,n]^2}{\longleftrightarrow}\{n\}\times[0,n] )\geq 1-\delta \}.
\end{equation*}
By results of \cite{Ke87}, there exists some $\delta_0$ such that $L(p,\delta_1)\asymp L(p,\delta_2)$ for $\delta_1,\delta_2<\delta_0$, so that we may refer unambiguously to $L(p)=L(p,\delta)$ with $\delta<\delta_0$. Define also
\begin{center}
$p_n=\underset{p}{\sup}\;\{L(p)>n\}$.
\end{center}

\noindent Let $\mathsf{A}(n,p)$ denote the event 
\begin{equation*}
\{\exists\text{ a $p$-open circuit $\mathcal{C}\in\mathsf{Ann}(n/4,n/2)$ such that }\mathcal{C}\underset{p-\text{open}}{\longleftrightarrow}\infty \}. 
\end{equation*}

\begin{lemm}
We have
\begin{equation*}
\underset{k\rightarrow\infty}{\lim}\;\underset{n>1}{\sup}\; \mathbb{P}(N>k)=0.
\end{equation*}
\end{lemm}

\noindent\textbf{Proof.  } Take $l>0$ and condition on the event $\mathsf{A}(n,p_{n/l})$: 
\begin{equation*}
\mathbb{P}(N_n>k) \leq \mathbb{P}(\mathsf{A}(n,p_{n/l})^c)+\mathbb{P}(N_n>k, \mathsf{A}(n,p_{n/l}))
\end{equation*}
By (2.4) in \cite{DSV09}, the first term is bounded from above by \\ $C_1\exp(-C_2\frac{n}{L(p_{n/l})})=C_1\exp(-C_2l)$ for some constants $C_1,C_2>0$. The BK inequality applied to $\{N_n=k\}$ gives
\begin{equation*}
\mathbb{P}(N_n>k, \mathsf{A}(n,p_{n/l}))<\mathbb{P}(S(n/2)\underset{p_{n/l}\text{-open}}{\longleftrightarrow} S(n))^k 
\end{equation*}
for the second term. Now 
\begin{equation*}\begin{split}
\mathbb{P}((S(n/2)\underset{p_{n/l}\text{-open}}{\longleftrightarrow} S(n))^c) & =\mathbb{P}(\exists\text{ a dual $p_{n/l}$-closed circuit in $\mathsf{Ann}(n/2,n)$} ) \\
& \geq \mathbb{P}(\{0\}\times[0,n/l]\underset{p_{n/l}\text{-closed}}{\longleftrightarrow} \{n\}\times[0,n/l])^4 \\
& \geq \delta(l)
\end{split}\end{equation*}
where the first inequality follows from the the FKG inequality \cite{W09} applied to gluing rectangle crossing events which imply the $p_{n/l}$-closed circuit. The existence of a lower bound $\delta(l)\searrow 0$ follows from the definition of $p_{n/l}$ and the Russo-Seymour-Welsh Theorem \cite{No08}. Thus 
\begin{equation*}
\mathbb{P}(N_n>k)\leq C_1\exp(-C_2l) + (1-\delta(l))^k.
\end{equation*}
We take $l$ large and $k\rightarrow \infty$. \\

\begin{lemm} Let $\omega$ be a configuration of open and closed edges and $e$ an edge in $\mathsf{Ann}(n/2,n)$. Define $\omega^+$ (resp. $\omega^-$) to be the configuration $\omega$ with the modification that $e$ have open (resp. closed) status. Suppose that $N_n(\omega^+)\neq N_n(\omega^-)$. Then there exists a dual circuit $\mathcal{C}\subset\mathsf{Ann}(n/2,n)$ about the origin, containing $e^*$ and having at most $N_n(\omega)$ open edges (we say $\mathcal{C}$ has at most $N_n$ \textit{defects}).\\
\end{lemm}

\noindent\textbf{Proof.  } By definition of $N_n$, we have $N_n(\omega^+)=N_n(\omega^-)+1$. Menger's theorem \cite{W09} applied to $\omega^-$ implies that there exists a closed dual circuit with $N_n(\omega^-)$ defects around the origin in $\mathsf{Ann}(n/2,n)$. Let $\omega^{--}$ be the configuration obtained by closing all edges on this circuit. Now $N_n(\omega^{--})=0$, so there exists at least one dual closed circuit $\mathcal{C}\subset\mathsf{Ann}(n/2,n)$ about the origin. Now define $\omega^{--+}$ to be the configuration obtained by closing $e$ in $\omega^{--}$. Suppose $N_n(\omega^{--+})=0$. Then the $N_n(\omega^-)$ edges which differentiate $\omega^{--+}$ and $\omega^{+}$ form a set which disconnects $\partial S(n/2)$ from $\partial S(n)$ in $\omega^{+}$, but by Menger's theorem this requires at least $N_n(\omega^+)>N_n(\omega^-)$ edges. Thus $N_n(\omega^{--+})>0$. Since $\omega^{--+}$ and $\omega^{--}$ differ only on $e$, we must have $e\in\mathcal{C}$ for some $\mathcal{C}$. Thus in $\omega^-$, $\mathcal{C}\setminus e$ contains at most $N_n(\omega^-)$ defects and the result follows.\\

\noindent Cover $T_\epsilon$ by $K_\epsilon$ balls $(B^j)_j$ of radius $\epsilon n$. For $j=0,1,\ldots K_\epsilon$, let $\omega_j$ be the configuration defined by 

\begin{equation}\begin{cases}
\omega_j(e) =\omega_{\epsilon}(e) \text{ for any edge } e\in \cup_{i=0}^j B_j\\
\omega_j(e) =\omega(e) \text{ otherwise}.
\end{cases}\end{equation}

\noindent It is clear that if $N$ changes after resampling $T_{\epsilon}$ (i.e. $N(\omega)\neq N(\omega_{\epsilon})$), then $N(\omega^{j-1})\neq N(\omega^j)$ for some $j$. \\

\begin{lemm}
Let $\bm{\upalpha}_4^K(p,q,s,n)$ denote the alternating four-arm event 
\begin{equation*} 
\left \{ \partial S(s) \xleftrightarrow[\text{having each at most $K$ defects }]{\text{4 paths alternating $p$-open and $q$-closed}} \partial S(s) + \partial S(n) \right \}
\end{equation*}

\noindent Then there exists a constant $C>0$ and a constant $C_l$ depending only on $l$ such that
\begin{gather*}
    \mathbb{P}(N(\omega^{j-1})\neq N(\omega^{j}),N(\omega^{j-1})\leq M, \mathsf{A}(n,p_{n/l})) \\
    \leq C_l (C\log \frac{1}{\epsilon})^{2M}\mathbb{P}(\bm{\upalpha}_4^0(p_c,p_c,\epsilon n,n))
\end{gather*}

\end{lemm}

\noindent \textbf{Proof.}\;\; Recall that the event $\mathsf{A}(n,p_{n/l})$ ensures the IPC adds only $p_{n/l}$-open edges after intersecting $S(n/2)$. On $\mathsf{A}(n,p_{n/l})$, the probability the infinite $p_{n/l}$-open path $y$ intersects $B^{j-1}$ and does not enter another ball after exiting $B^{j-1}$ is $O(1/K_{\epsilon})$, so we may assume this and obtain the statement up to an additive constant of order $\epsilon$. Let $S$ be a box of side length $n/2$ centered about $B^{j-1}$, $s_1,s_2$ be the disjoint components of $S\setminus y$, and $S^+,S^-=s_1\setminus \text{int } B^{j-1}, s_2\setminus \text{int } B^{j-1}$. Then, on $\{N(\omega)\leq M\}$, there can not exist $M$ disjoint $p_c$-open paths $\subset \text{cl } S^+$ connecting the disjoint segments $y_1,y_2$ of $y\setminus\text{cl } B^{j-1}$. To see this, apply Lemma 4.3 to each edge $e\in B^{j-1}$ to obtain the inclusion

\begin{multline*}
\quad\quad\quad\quad\quad\quad\quad 
\big \{ N(\omega)\leq M \big \} \cap \big \{N(\omega^{j-1})\neq N(\omega^{j}) \big \}  \\
\subset \\ 
\Bigg \{ \exists\text{ a closed  circuit }\mathcal{C}\text{ intersecting }B^{j-1} \text{ in }\mathsf{Ann}(n/2,n)^* 
\text{ with $\leq M$ defects} \Bigg \}
\end{multline*}
 
\noindent so that $\mathcal{C}\setminus B^{j-1}$ contains at most $M-1$ defects. Now if there did exist $M$ such paths, these paths would include an edge of $\mathcal{C}\setminus\{\text{defects}\}$, which would thus be $p_c$-open. Given $\mathsf{A}(n,p_{n/l})$, this edge would have to be open, contradicting that $\mathcal{C}$ is closed. By Menger's Theorem, there exists a dual $p_c$-closed path with at most $M$ defects joining $\partial B^{j-1}$ to $\partial S^+\setminus\cup_i y_i$. The same argument applied to $S^-$ gives a dual $p_c$-closed path joining $\partial B^{j-1}$ to the other side of $S$. These observations imply the event $\bm{\upalpha}_4^M(p_{n/l},p_c,\epsilon n,n/4)$.\\
\begin{figure}
\vspace{-150pt}
\hspace{-40pt}
\includegraphics[scale=.8]{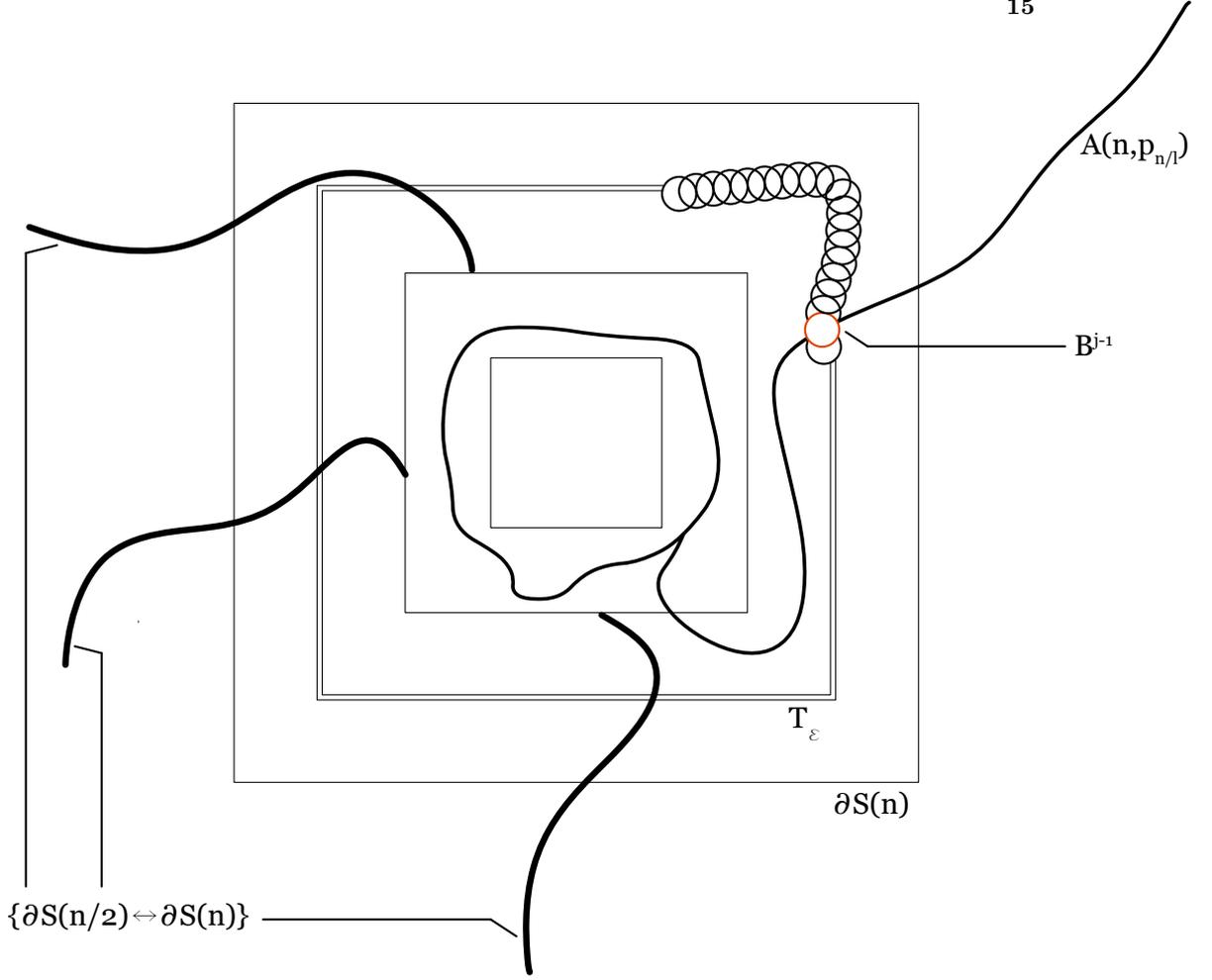}
\vspace{-150pt}
\caption{The setup for Lemma 4.4. The infinite path from the event $\mathsf{A}(n,p_{n/l})$ passes through the covering of $T_{\epsilon}$, and disjoint arms for $N(\omega)$ cross $\mathsf{Ann}(n/2,n)$. The dual $p_c$-closed circuit enclosing $S(n/2)$ leading to the event $\bm{\upalpha}_4^M(p_{n/l},p_c,\epsilon n,n/4)$ is omitted.}
\end{figure}

\noindent Next, following the argument in Proposition 18 in \cite{No08}, the condition of having $\leq M$ defects in $\bm{\upalpha}_4^M(p_{n/l},p_c,\epsilon n, n/4)$ can be removed at the cost of logarithmic factor:
\begin{equation*}
\mathbb{P}(\bm{\upalpha}_4^M(p_{n/l},p_c,\epsilon n,n))\leq C^{2M}(1+\log(n/\text{rad }B))^{2M}\mathbb{P}(\bm{\upalpha}_4^0(p_{n/l},p_c,\epsilon n, n))
\end{equation*}
This follows from that Proposition by modifying the induction hypothesis on $M$ to establish the upper bound
\begin{equation}
\mathbb{P}(\bm{\upbeta}_4^M(p_{n/l},p_c,\epsilon n, n))\leq C_M \mathbb{P}(\bm{\upbeta}_4^0(p_{n/l},p_c,\epsilon n, n)
\end{equation}
such that the factor $C_M$ be taken in the form $C^M$. Here the event $\bm{\upbeta}_4^M(p,q,s,n)$ is defined similarly to $\bm{\upalpha}_4^M(p,q,s,n)$ except that the total number of defects is $M$ rather than $2M$. \\

\begin{figure}
\vspace{-200pt}
\hspace{-50pt}
\includegraphics[scale=0.79]{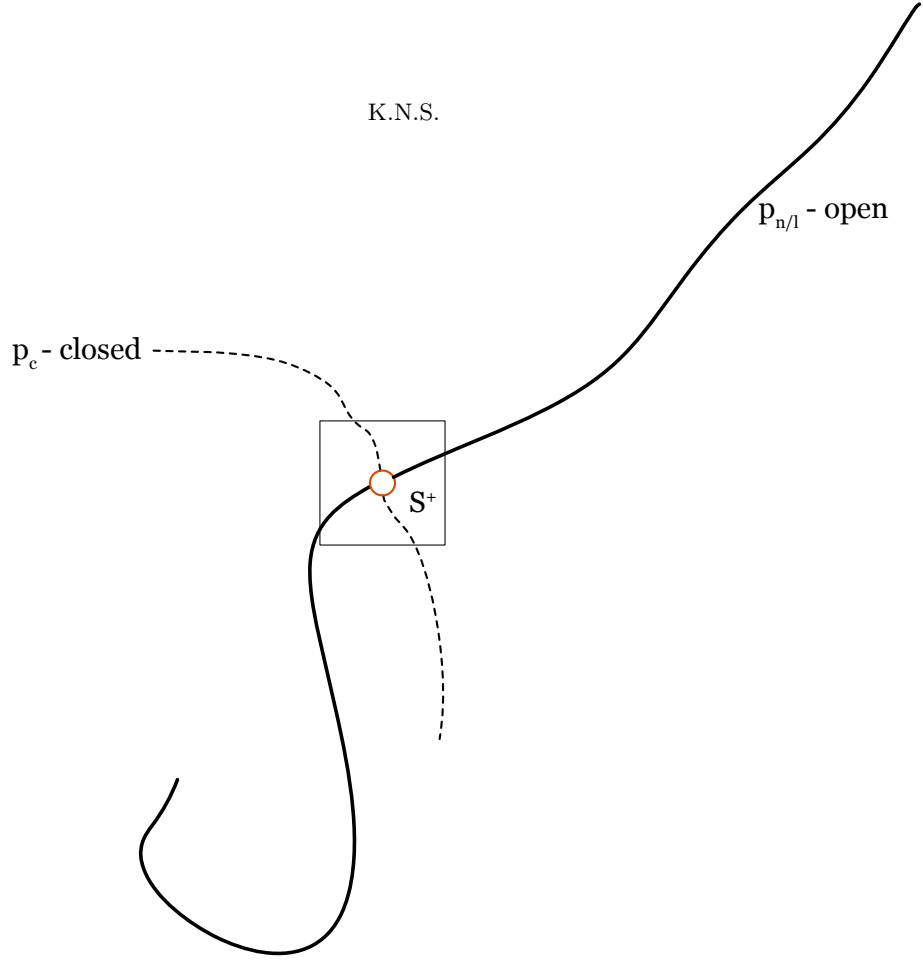}
\vspace{-150pt}
\caption{The alternating four-arm event for $B^{j-1}$ (\textit{ex.} defects).}
\end{figure}

\noindent Now by inclusion and the definition of correlation length,
\begin{equation*}\begin{split}
\mathbb{P}(\bm{\upalpha}_4(p_{n/l},p_c,\epsilon n,n/4)) & \leq \mathbb{P}(\bm{\upalpha}_4(p_{n/l},p_c,\epsilon n,n/l)) \asymp \mathbb{P}(\bm{\upalpha}_4(p_c,p_c,\epsilon n, n/l)) \\
& \leq \text{const} \frac{\mathbb{P}(\bm{\upalpha}_4(p_c,p_c,\epsilon n, n/4))}{\mathbb{P}(\bm{\upalpha}_4(p_c,p_c,n/l,n/4))} \\
& \leq C_l \mathbb{P}(\bm{\upalpha}_4(p_c,p_c,\epsilon n, n/4)) 
\end{split}\end{equation*}
where the two last inequalities follow from quasi-multiplicity of four-arm probabilities and the RSW Theorem applied to $\bm{\upalpha}_4(p_c,p_c,n/l,n/4)$ \cite{No08}. \\

\begin{lemm}
There exists a constant $C_l$ depending only on $l$ such that
\begin{equation*}
\mathbb{P}(\mathsf{A}(n,p_{n/l})(\omega)\cap \mathsf{A}(n,p_{n/l})^c(\omega^j))\leq C_l \mathbb{P}(\bm{\upalpha}_4(p_c,p_c,\epsilon n,n/4))
\end{equation*}
\end{lemm} 

\noindent \textbf{Proof.}\;\; The event $\left \{\mathsf{A}(n,p_{n/l})(\omega)\cap \mathsf{A}(n,p_{n/l})^c(\omega^j) \right \}$ implies $B^j$ is connected by 4 $p_{n/l}$-open arms to $\partial S(n)$. By the definition of correlation length, 
\begin{equation*}
\bm{\upalpha}_4(p_{n/l},p_{n/l},\epsilon n,n/4)\leq \bm{\upalpha}_4(p_{n/l},p_{n/l},\epsilon n,n/l) \asymp \mathbb{P}(\bm{\upalpha}_4(p_c,p_c,\epsilon n,n/l))
\end{equation*}
and the statement follows from quasi-multiplicity and the RSW Theorem. \\

\noindent \textbf{Proof of Theorem 4.1.}\;\; Conditioning on $\mathsf{A}(n,p_{n/l})$,
\begin{equation*}\begin{split}
& \mathbb{P}(\exists\; j:\; N(\omega^{j-1})\neq N(\omega^j) )  \leq \mathbb{P}(\mathsf{A}(n,p_{n/l})^c) + \mathbb{P}(N(\omega)>M) \\
+ &  \mathbb{P}(\exists\; j:\; N(\omega^{j-1})\neq N(\omega^j), \mathsf{A}(n,p_{n/l})(\omega),N(\omega)<M).
\end{split}\end{equation*}

\noindent The third term is bounded from above by
\begin{equation*}\begin{split}
& \mathbb{P}(\mathsf{A}(n,p_{n/l})(\omega)\cap \cup_j \mathsf{A}(n,p_{n/l})^c(\omega^j)) \\
+ & \mathbb{P}(\exists\; j: N(\omega^{j-1}) \neq N(\omega^{j}), \mathsf{A}(n,p_{n/l})(\omega^s)\;\forall\; s, N(\omega)\leq M)
\end{split}\end{equation*}
and by Lemma 4.5,
\begin{equation*}\begin{split}
& \mathbb{P}(\mathsf{A}(n,p_{n/l})(\omega)\cap \cup_j \mathsf{A}(n,p_{n/l})^c(\omega^j))\\
\leq & \sum_{j=1}^{K_{\epsilon}} \mathbb{P}(\mathsf{A}(n,p_{n/l})(\omega)\cap \mathsf{A}(n,p_{n/l})^c(\omega^j)) \\
\leq & C_l K_{\epsilon}\mathbb{P}(\bm{\upalpha}_4(p_c,p_c,\epsilon n, n/4)).
\end{split}\end{equation*} By Lemma 4.3,
\begin{equation*}\begin{split}
& \mathbb{P}(\exists\; j: N(\omega^{j-1}) \neq N(\omega^{j}), \mathsf{A}(n,p_{n/l})(\omega^s)\;\forall\; s, N(\omega)\leq M) \\
\leq & \sum_{j=1}^{K_{\epsilon}} \mathbb{P}(N(\omega^{j-1}) \neq N(\omega^{j}), \mathsf{A}(n,p_{n/l})(\omega^{j-1}), N(\omega)\leq M) \\
\leq & C_l K_{\epsilon} (C\log \frac{1}{\epsilon})^{2M} \bm{\upalpha}_4(p_c,p_c,\epsilon n,n/4).
\end{split}\end{equation*}
By Lemma 4.2, 
\begin{equation*}
\underset{M\rightarrow \infty}{\limsup}\; \mathbb{P}(\mathsf{A}(n,p_{n/l})^c) + \mathbb{P}(N(\omega)>M)=0
\end{equation*}
and to conclude (4) we recall that $K_{\epsilon}$ may be taken $\leq \text{const}/\epsilon$ by regularity of the resampled region $T_{\epsilon}$, and that $\bm{\upalpha}_4(p_c,p_c,\epsilon n, n/4)\leq \text{const }\epsilon^{\alpha_4}$, where the critical four-arm exponent $\alpha_4$ is $>1$ \cite{SW01}. $\blacksquare$  
\\ \\

\noindent\textbf{Sensitivity and Factorization}\;\; It should be possible to follow the above approach of using near-critical percolation arguments to obtain a mesh-independent version of discrete splicing, and thus factorization of invasion on planar annuli. This would provide one description of the scaling limit of invasion percolation, i.e., by indexing the process by annuli crossed by a specified number of disjoint paths, and complements the approach indicated in \cite{GPS09} for describing scaling limits of near-critical and dynamical percolation. One natural direction is to ask similar questions for the minimal spanning tree, which may inherit sensitivity properties from invasion as a subgraph, but also carries translational invariance, which would make it a good candidate for a two-dimensional black noise. \\ \\ 

\noindent\textit{Acknowledgements}. I would like to thank Michael Damron for many helpful discussions about invasion percolation and for supervising this work, Pierre Bertin for technical explanations during a reading project on black noise, and Ya. Sinai for valuable intuition in his random processes course.

\medskip\medskip\medskip\medskip\medskip\medskip\medskip\medskip\medskip\medskip\medskip\medskip\medskip\medskip\medskip\medskip\medskip\medskip\medskip\medskip\medskip\medskip\medskip\medskip\medskip\medskip\medskip\medskip\medskip\medskip\medskip\medskip\medskip\medskip\medskip\medskip\medskip\medskip\medskip\medskip\medskip\medskip\medskip\medskip\medskip\medskip\medskip\medskip\medskip\medskip\medskip\medskip\medskip\medskip\medskip\medskip\medskip\medskip\medskip\medskip\medskip\medskip\medskip\medskip\medskip\medskip\medskip\medskip\medskip\medskip\medskip\medskip\medskip\medskip\medskip\medskip\medskip\medskip\medskip\medskip\medskip\medskip\medskip\medskip


\bibliography{bibtemplate}
\bibliographystyle{smfalpha} 

\end{document}